\documentclass[a4paper,12pt]{amsart}
\usepackage [latin1]{inputenc}
\usepackage[all]{xy}
\usepackage{amsmath}
\usepackage{amsfonts}
\usepackage{amssymb}
\usepackage{graphics}
\usepackage{xy,amsmath,amssymb}
\usepackage[dvips]{graphicx}
\usepackage{graphicx}
\usepackage[initials, non-compressed-cites]{amsrefs}
\newtheorem{theorem}{Theorem}[section]
\newtheorem*{theorem*}{Theorem}
\newtheorem{lemma}[theorem]{Lemma}
\newtheorem{prop}[theorem]{Proposition}
\newtheorem{cor}[theorem]{Corollary}
\newtheorem{ex}{Example}[section]
\newtheorem{remark}{Remark}[section]
\newtheorem{defin}[theorem]{Definition}

\def\be{\begin{equation*}}
\def\ee{\end{equation*}}
\def\bt{\begin{theorem}}
\def\et{\end{theorem}}
\def\bp{\begin{prop}}
\def\ep{\end{prop}}
\def\bl{\begin{lemma}}
\def\el{\end{lemma}}
\def\bc{\begin{cor}}
\def\ec{\end{cor}}
\def\br{\begin{oss}\rm}
\def\er{\end{oss}}
\def\br{\begin{remark}}
\def\er{\end{remark}}
\def\bex{\begin{ex}\rm}
\def\eex{\end{ex}}
\def\bd{\begin{defin}}
\def\ed{\end{defin}}
\def\demo{\par\noindent{\bf Proof.\,\,}}
\def\bit{\begin{itemize}} \def\eit{\end{itemize}}
\def\ben{\begin{enumerate}} \def\een{\end{enumerate}}
\def\smi{\setminus}
\def\enddemo{\ $\Box$\par\vskip.6truecm}
\def\R{{\mathbb R}}   \def\a {\alpha} \def\b {\beta}\def\g{\gamma}
\def\N{{\mathbb N}}     
\def\e{\varepsilon}
\def\C{{\mathbb C}}
     \def\l{\lambda}
                 \def\o{\omega}\def\p{\partial}
\def\Z{{\mathbb Z}}       \def\s{\sigma} \def\z{\zeta}
 \def\tms{\times} 
\def\sbs{\subset}    
\def\ES{\varnothing}
 \def\oli{\overline}

  \def\oli{\overline}
\def\G{\Gamma}
\def\nin{\noindent}
\def\til{\tilde} \def\Til{\widetilde}
\def\O{\Omega}  \def\beq{\begin{eqnarray}}
\def\eeq{\end{eqnarray}}
\def\oli{\overline}
\def\IN{\infty}

\def\({\rm(}
\def\){\rm)}

\def\gfra{\mathfrak{g}}

\newcommand{\lra}{\longrightarrow}

\newcommand{\lcal}{\mathcal{L}}
\newcommand{\mcal}{\mathcal{M}}

\newcommand{\ocal}{\mathcal{O}}
\newcommand{\hcal}{\mathcal{H}}

\newcommand{\zcal}{\mathcal{Z}}
\newcommand{\fcal}{\mathcal{F}}
\newcommand{\gcal}{\mathcal{G}}
\newcommand{\ucal}{\mathcal{U}}

\newcommand{\ccal}{\mathcal{C}}

\newcommand{\ecal}{\mathcal{E}}
\newcommand{\acal}{\mathcal{A}}

\newcommand{\tcal}{\mathcal{T}}
\newcommand{\scal}{\mathcal{S}}

\newcommand{\rcal}{\mathcal{R}}

\def\benu{\begin{enumerate}} \def\eenu{\end{enumerate}}
\def\beqn{\begin{eqnarray*}}  \def\eeqn{\end{eqnarray*}}
\def\beqnn{\begin{eqnarray*}}  \def\eeqnn{\end{eqnarray*}}

\def\rdim{{\rm dim}_{\,_\R}}

\def\Gr{\mathrm{Gr}}
\def\rk{\mathrm{rk}}

\def\({\mathrm{(}} \def\){\mathrm{)}}
\title[Oka principle for Levi flat manifolds]{Oka principle for Levi flat manifolds}

\author[Samuele Mongodi and Giuseppe Tomassini ]{Samuele Mongodi and Giuseppe Tomassini }
              \address{Politecnico di Milano, Dipartimento di Matematica, Via Bonardi, 9 -- I-20133 Milano, Italy}
              \email{samuele.mongodi@polimi.it
              }
              \address{Scuola Normale Superiore, Piazza dei Cavalieri, 7 -- I-56126 Pisa, Italy
              }
              \email{g.tomassini@sns.it
              }

 \date{\today}


             \dedicatory{In memory of Paolo, colleague and friend}

\begin{document}

\maketitle

 \section{Introduction}

The name of \emph{Oka principle}, or \emph{Oka-Grauert principle}, is traditionally used to refer to the holomorphic incarnation of the homotopy principle: on a Stein space, every problem that can be solved in the continuous category, can be solved in the holomorphic category as well; this line of thought originated from a paper by Oka \cite{oka}, where he shows that a topologically trivial line bundle on a domain of holomorphy is holomorphically trivial. The underlying idea was further explored in Grauert's work on the classification of holomorphic fiber bundles \cites{G1,G2,G3}; subsequently, inspired by Gromov's seminal paper \cite{gro}, Forstneric and others developed the Oka principle into a well formed and exhaustive theory (see \cite{forst}).

\medskip
In this note, we begin the study of the same kind of questions on a Levi-flat manifold; more precisely, we try to obtain a classification of CR-bundles on a semiholomorphic foliation of type $(n,1)$. Our investigation should only be considered a preliminary exploration, as it deals only with some particular cases, either in terms of regularity or bidegree of the bundle, and partial results.

In order to make our intent clearer, we anticipate some of the results and notions presented in the paper. We refer to Sections 2 and 3 for the precise definitions of the objects involved.

Given a (smooth) semiholomorphic foliation $X$ of type $\(n,d\)$ let
$${\sf Vect}_{\,\sf top}^{(m,l)}(X):=\Big\{{\rm topological\,\, vector\,\, bundles\,\, of\,\, birank\,\,} (m,l)\Big\}\;,$$
$${\sf Vect}_{\,\sf cr}^{(m,l)}(X):=\Big\{{\rm smooth\,\, CR\,\, vector\,\, bundles\,\, of\,\, birank\,\,} (m,l)\Big\}\;.$$
Elements of ${\sf Vect}_{\,\sf cr}^{(m,l)}(X)$ are called CR {\em vectors bundles}.

If $X$ is real analytic let
$${\sf Vect}_{\,{\sf cr},{\sf an}}^{(m,l)}(X):=\Big\{{\rm real\,\,analytic\,\, CR\,\, vector\,\, bundles\,\, of\,\, birank\,\,} (m,l)\Big\}\;.$$
The CR analogous of the Oka-Grauert principle can be formulated in the following way:

\begin{quote}for strongly 1-complete semiholomorphic foliations $X$ of type $\(n,1\)$ the natural map
\be
{\sf\epsilon}_{\!_X}:{\sf Vect}_{\,\sf cr}^{(m,l)}(X)\longrightarrow
{\sf Vect}_{\,\sf top}^{(m,l)}(X)
\ee
is a bijection.
\end{quote}

Precise definitions can be found in section 2. Let us immediately give a counterexample to the real analytic case of the result.

\bex
Let $X=\C\tms\R$, $\mathcal S^\o$ be the sheaf of germs of real analytic CR functions in $X$. Then, since $H^1(X,\Z)=H^2(X,\Z)=0$ we have
$$
H^1(X,\mathcal S^\o)\simeq H^1(X,\mathcal S^{\o\,\ast})
$$
i.e. the real analytic CR bundles of birank $\(1,0\)$ are (topologically trivial and) parametrized by the cohomology group $H^1(X,\mathcal S^\o)$ which is $\neq 0$ (see \cite{annac}). Thus the Oka-Grauert principle fails in ${\sf Vect}_{\,{\sf cr},{\sf an}}^{(m,l)}(X)$.
\eex

As in the holomorphic category, the validity of the Oka-Grauert principle can be rephrased in terms or homotopy of CR maps.

Indeed, every CR vector bundle ${\sf E}$ of birank $(m,l)$ embeds, as a topological bundle in the trivial bundle $X\tms(\C^M\tms\R^L)$ for some $M, L\in\Z$. Then, ${\sf E}\simeq f^\ast {\sf U}^{M,L}_{m,l}$, where $f:X\to G^{M,L}_{m,l}$ sends $x\in X$ to $E_x$, $G^{M,L}_{m,l}$ and ${\sf U}^{M,L}_{m,l}$ being appropriate CR versions of the Grassmannian and the universal bundle - we refer to section 3 for the definitions.  Therefore surjectivity of ${\sf\epsilon}_{\!_X}$ will be be a consequence of the following assertion:
\begin{quote}
(i) every continuous map $f:X\to G^{M,L}_{m,l}$ is homotopic to a CR map.
\end{quote}

As for injectivity let ${\sf E},{\sf E}'\in{\sf Vect}_{\,\sf cr}^{(m,l)}(X)$, $\{h_{ij}\}$, $\{h'_{ij}\}$ their cocycles (associated to a trivializing open covering $\{U_{i}\}$). Let $Z_i:=U_i\tms{\sf G}_{m,l}$ (${\sf G}_{m,l}$ is an appropriate group of matrices, see section 2) and $Z:=\amalg_{i}Z_i/\!\!\sim$ where $\sim$ identifies $(x,v)$ to $(x,v')$, $x\in U_i\cap U_j$, where $v'=h'_{ij}(x)vh_{ji}(x)$.

Then, $Z\to X$ is a CR bundle with fiber ${\sf G}_{m,l}$ and as it easily seen, topological (respectively CR) isomorphisms ${\sf E}\simeq {\sf E}'$ correspond to continuous (respectively CR) sections of $Z\to X$. Therefore the sentence ${\sf\epsilon}_{\!_X}$ is injective  will be a consequence of the following assertion:
\begin{quote}
(ii) every continuous section $X\to Z$ is homotopic to a CR section.
\end{quote}

However, a number of technical difficulties arise, when trying to adapt this line of proof to the CR case, so our approach steers, in part, toward cohomological methods.

We present the main definitions in section 2, recalling the concepts of semiholomorphic foliation and CR-bundles from our previous work \cite{mt}; in section 3, we define the CR Grassmannian and we show that it has a natural structure of semiholomorphic foliation. Section 4 recalls the concept of 1-completeness from \cite{mt} and presents some vanishing results for cohomology.

Sections 5 and 6 contain the results we obtained on the Oka-Grauert principle. In section 5, we examine the case of real analytic vector bundles: as we pointed out, we cannot hope for an equivalence between topological and CR real analytic classifications, nonetheless the real analytic structure proves to be of help in proving that the topological and the CR smooth classifications are equivalent. We have the following result (see Theorems \ref{teo3}, \ref{teo31} and \ref{teo32}):
\begin{theorem*}Let $X$ be a strongly $1$-complete real analytic semiholomorphic foliation of type $(n,1)$. Let ${\sf E}\to X$, ${\sf F}\to X$ be topologically equivalent real analytic $CR$ bundles of birank $(m,l)$. Then:
\begin{enumerate}
\item if ${\sf E}\to X$ is topologically trivial, then it is (smoothly) CR trivial;
\item if $l=0$, ${\sf E}\to X$ and ${\sf F}\to X$ are (smoothly) CR equivalent;
\item if $l=1$ and $H^1(X,\mathcal{T})=H^1(X,\mathbb{Z})=0$, then ${\sf E}\to X$ and ${\sf F}\to X$ are (smoothly) CR equivalent.
\end{enumerate}
\end{theorem*}
The sheaf $\mathcal{T}$ is the sheaf of germs of smooth functions which are constant on the leaves of $X$.

In section 6, we tackle the problem for smooth CR bundles: we manage to give a complete answer for $(1,0)$ bundles (Theorem \ref{teo1OG}), reducing the problem for $(1,1)$ to a cohomological property (resembling of the approach to the classical problem presented in \cite{Cart}). We also indicate a possible route to deal with the general case of $(m,l)$ bundles, by reducing it to the corresponding problem for $(m,0)$ and $(0,l)$ bundles.

\section{Preliminaries on semiholomorphic foliations}

We will briefly introduce the objects that we are going to use and study in this note; we refer also to \cites{mt,mt1} for a broader discussion of semiholomorphic foliations.

\subsection{Main definitions}We recall that a {\it semiholomorphic foliation of type $(n,d)$} is a (connected) smooth foliation $X$ whose local models are subdomains $U_j=V_j\tms B_j$ of $\C^n\tms\R^d$ and whose local changes of coordinates $(z_k,t_k)\mapsto(z_j,t_j)$ are of the form
\begin{equation}\label{LT}
\begin{cases}
z_j=f_{jk}(z_k,t_k) \>\> &\\
t_j=g_{jk}(t_k),\>\> &
\end{cases}
\end{equation}
where $f_{jk}$, $g_{jk}$ are smooth and $f_{jk}$ is holomorphic with respect to $z_k$. If we replace $\R^d$ by $\C^d$ and we suppose that $f$ and $h$ are holomorphic we get the notion of {\it holomorphic foliation of codimension} $d$.

Local coordinates $z^1_j,\ldots,z^n_j,t^1_j,\ldots,t^d_j$ satisfying (\ref{LT}) are called {\it local distinguished coordinates}.

A closed subset $Y\sbs X$ is said to a {\it subsemiholomorphic foliation} of $X$ if for every point $p\in Y$ there exist local distinguished coordinates $z^1,\ldots,z^n,t^1,\ldots,t^d$ on a neighborhood $U$ of $y$ in $X$ such that
$$
U\cap Y=\Big\{x\in U:z_{m+1}=\ldots=z_n=t_{s+1}=\ldots=t_d=0\Big\}.
$$
It follows that $Y$ is a semiholomorphic foliation of type $(m,s)$.
\bex $\C\tms\R$ is not a subsemiholomorphic foliation of $\C^2$ considered as a foliation of type $(2,0)$ but it is if  $\C^2$ is considered as a foliation of type $(1,2)$.\eex

We denote $\scal=\mathcal S_X$ the sheaf of germs of smooth CR functions in $X$. If $X$ is real analytic we define $\mathcal S^\o$, the sheaf of germs of real analytic CR functions in $X$. For every open set $U$ in $X$ we consider on $\mathcal S_X(U)$  the Fr\'echet topology induced by $C^\IN(U)$. Clearly, $\mathcal S(U)$ is closed subset. If $X$, $Y$ are semiholomorphic foliations of types $(m,s)$, $(m',s')$ a map $f:X\to Y$ is said to be CR if $f_\ast\mathcal S_Y\sbs \mathcal S_X$. We denote $\ccal\rcal(X,Y)$ the set of all CR maps $X\to Y$.

The sheaf of germs of smooth functions which are locally constant on the leaves is denoted by $\tcal=\mathcal T_X$. If $X$ is real analytic, $\mathcal T^\o\sbs\mathcal T$ is  the subset of those germs which are real analytic.

\subsection{Complexification} A real analytic semiholomorphic foliation of type $(n,d)$ can be {\em complexified}, essentially in a unique way: there exists a holomorphic foliation $\Til X$ of type $(n,d)$ with a closed real analytic CR embedding $X\hookrightarrow\Til X$ (cfr. \cite{rea}*{Theorem 5.1}). In particular, $X$ is  a Levi flat submanifold of $\Til X$

In order to construct $\Til X$, we consider a covering by distinguished domains $\{U_j=V_j\tms B_j\}$ and we complexify each $B_j$ in such a way to obtain domains $\Til U_j$ in $\C^n\tms\C^d$. The domains $\Til U_j$ are patched together by the local change of coordinates
\be
\begin{cases}
z_j=\Til f_{jk}(z_k,\tau_k) \>\> &\\
\tau_j=\Til g_{jk}(\tau_k)\>\> &
\end{cases}
\ee
obtained complexifying the (vector) variable $t_k$ by $\tau_k=t_k+i\theta_k$ in the real analytic functions $f_{jk}$ and $g_{jk}$ (cfr. \eqref{LT}).

In the sequel, we will call complexification of $X$ every neighborhood of $X$ in $\Til X$.

\subsection{CR-bundles}\label{crbu}
Let ${\sf G}_{m,l}$ be the group of matrices
$$
\left( \begin{array}{cc}
	A&B\\
	0&C
\end{array} \right)
$$
where $A\in {\rm GL(m,\C)}$, $B\in {\rm GL(m,l,\C)}$ (the set of $m\tms l$ matrices with complex entries) and $C\in {\rm GL(l,\R)}$.  We set ${\sf G_{m,0}}={\rm GL(m,\C)}$, ${\sf G_{0,l}}={\rm GL(l,\R)}$.

To each matrix $M\in \sf G_{m,l}$ we associate the linear transformation $\C^m\tms\R^l\rightarrow\C^m\tms\R^l$ given by
 $$
 (\xi,u)\mapsto (A\xi+Bu,Cu),
 $$
$(\xi,u)\in\C^m\tms\R^l$.

Let $X$ be a semiholomorphic foliation of type $(n,d).$  A CR {\it bundle of birank (m,l)} is a vector bundle $\pi:{\sf E}\to X$ whose cocycle with respect to a countable, trivializing distinguished covering $\{U_j\}$ is a set of smooth CR map
$\g_{jk}:U_j\cap U_k\to\sf G_{m,l}$
\begin{equation}\label{CO}
\g_{jk}=\left( \begin{array}{cc}
	A_{jk}&B_{jk}\\
	0&C_{jk}
\end{array} \right)
\end{equation}
where $C_{jk}=C_{jk}(t)$ is a matrix with smooth entries and $A_{jk}$, $B_{jk}$ are matrices with smooth CR entries. Thus $E$ is foliated by complex leaves of dimension $n+m$ and real codimension $d+l$.

We denote $\ecal=\ecal_X$ the sheaf of germs of smooth CR sections of ${\sf E}$ and we set $\ecal^\o=\ecal\cap C^\o$ whenever $X$ and ${\sf E}$ are real analytic. If $U$ is an open subset of $X$ the topology on $\ecal(U)$, the space of CR sections $s:U\to {\sf E}$ is defined as follows. We fix a trivializing distinguished covering $\{U_j\}_{\in\N}$ of $X$ and let $\{\g_{jk}\}$ be the corresponding cocycle. Then $\ecal(U_j\cap U)\simeq\ccal\rcal(U_j\cap U,\C^m\tms\R^l)$ therefore $\ecal(U)$ identifies to
$$
 \Big\{(f_j)\in\prod\limits_{j\in\N}\ccal\rcal(U_j\cap U,\C^m\tms\R^l                                  ):{f_j}_{|U_j\cap U_k\cap U}=\g_{jk} {f_k}_{|U_j\cap U_k\cap U}\Big\}
$$
which is (a closed subspace of) a Fr\'echet space.
Let $X$ be a semiholomorphic foliation of type $(n,d)$. Then
\bit
\item the tangent bundle ${\sf T}X$ of $X$ is a CR-bundle of birank $(n,d).$
\item the bundle ${\sf T}_\lcal={\sf T}_\lcal^{1,0}X$ of the holomorphic tangent vectors to the leaves of $X$ is a CR-bundle of type $(n,0).$
\item the transverse bundle ${\sf N}_{\rm tr}$ (to the leaves of $X$) is a CR-bundle of type $(0,d).$
\eit
Observe that If $X$ is embedded in a complex manifold $Z$, its transverse ${\sf T}Z/{\sf T}X$ is not a CR bundle in general.

If ${\sf E}$ is a real analytic CR bundle of type $(m,l)$, given by the cocycle \eqref{CO}, then complexifying $\g_{jk}$ by
$$
\Til\g_{jk}=\left( \begin{array}{cc}
	\Til A_{jk}&\Til B_{jk}\\
	0&\Til C_{jk}
\end{array} \right)
$$
and the fibre $\C^m\tms \R^l$ by $\C^m\tms \C^l$ we obtain a holomorphic vector bundle $\Til {\sf E}\to \Til X$ of rank $(m,l)$. We will call $\Til {\sf E}$ the {\em complexification} of $E$. A section $s$ of $\Til {\sf E}$ is locally given by a couple $(f_i,g_i)$ where $f_i:U_i\to\C^m$, $g_i:U_i\to\C^l$ such that
$$
\begin{cases}
f_i=\Til A_{ij}f_j+\Til B_{ij}g_j \>\> &\\
g_i=\Til C_{ij}g_j\>\> &
\end{cases}.
$$
If $s=\{(f_i,g_i)\}_i$ is a section of $\Til {\sf E}$ then $s_{\vert X}:=\{({f_i}_{\vert X},{{\sf Re}\,g_i}_{\vert X})\}_i$  is a section of ${\sf E}$.

\section{The CR Grassmannian}\label{CR-Gras}
Let $\Gr_k(\C^d)$ (resp. $\Gr_k(\R^d)$) denote the grassmannian of the $k$-dimensional complex (real) subspaces of $\C^d$ (resp. $\R^d$).

Given two pairs $(m,l)$,$(M,L)$ of integers with $0\le m\le M, 0\le l\le L$, let \\$\Gr_{m,l}(\C^M\times\R^L)$ be the set of the real (linear) subspaces $V\subseteq \C^M\times\R^L$ such that
$$
\rdim V=2m+l,\>\>\>\dim_\C(V\cap (\C^M\times\{0\}))=m.
$$
By definition, $\Gr_{m,l}(\C^M\times\R^L)$, also denoted $\Gr^{M,L}_{m,l}$, is the \emph{Grassmannian variety} of the CR subspaces of $\C^M\times\R^L$ of type $(m,l)$. We want to prove that $\Gr_{m,l}(\C^M\times\R^L)$ is a semiholomorphic foliation.

We consider every $V\in \Gr_{m,l}(\C^M\times\R^L)$ as a subspace of $\C^M\times \C^L=\C^{M+L}$ and we denote by $V^\C=V\oplus i V$; $V^\C$ is an $(m+l)$-dimensional complex subspace of  so we have a map
$$
\Gr_{m,l}(\C^M\times\R^L)\to \Gr_{m+l}(\C^{M+L})
$$
sending $V$ in $V^{\C}$.
\subsection{Chart description}
The usual charts for the Grassmannian are given as follows: fix a basis $\{e_1,\ldots, e_{M+L}\}$ for $\C^{M+L}$ (so that the first $M$ vectors span $\C^M$ and the second $L$ span $\C^L$), for every $I\subseteq\{1, \ldots, M+L\}$ with $|I|=m+l$, let
$$V_{I}^0=\mathrm{span}_{\C}\{e_j,\ j\not\in I\}$$
and consider
$$
U_I=\{\Lambda\in \Gr_{m+l}(\C^{M+L})\ :\ \Lambda\cap V_{I}^0=\{0\}\}\;.
$$
$\Lambda\in U_I$ is then represented by the unique $(M+L)\times(m+l)$ matrix $A_\Lambda$, with the $(m+l)\times(m+l)$ minor corresponding to the multi-index $I$ is the identity and such that the image of $A_\Lambda:\C^{m+l}\to\C^{M+L}$ is $\Lambda$.

The other $(m+l)(M+L-m-l)$ entries of $A_\Lambda$ are the local coordinates for $U_I$.

We cut out a real analytic subvariety $\zcal$ of $U_I$ with the following equations:
$$a_{jk}=0\qquad  j=M+1,\ldots, M+L,\quad k=1,\ldots, m$$
$$a_{jk}=\oli{a}_{jk}\qquad j=M+1,\ldots, M+L,\quad k=m+1,\ldots, m+l\;.$$

We note that, if $I$ contains less than $m$ indexes in $[1,M]$, then no element of $U_I$ can satisfy the first set of equations, because they would contradict the definition of $U_I$.

\bigskip

To calculate the dimension of $\zcal$, let us note that, as the $(m+l)\times(m+l)$ minor defined by $I$ is the identity, then some of the equations written above are redundant.

We stratify $\zcal$ as follows; for $A\in \zcal$, let $A_s$ be the submatrix formed by the last $l$ columns and $L$ rows, then
$$
\zcal^r=\{A\in X\ :\ \rk( A_s)=r\}\;.
$$
Obviously, the rank of $A_s$ can never be more than $l$, as the rank of the submatrix given by the first $m$ columns and $M$ rows is at least $m$.

We note that $\zcal^r$ does not intersect $U_I$ if $I$ contains less than $r$ indices in $[M+1, M+L]$.

If we are at a point $x\in \zcal^r$, we will have $mr$ equations of the first set already satisfied and $lr$ equations already satisfied in the second set.

Therefore,
$$
\dim_{\R} \zcal^r=2m(M-m)+2l(M-m)+l(L-l)\;.
$$

The factors $2$ take into account the fact that some of the equations describe complex manifolds.

\bigskip

So, $\zcal^l$ is the top dimensional stratum of $\zcal$, i.e. the regular points of the irreducible components of top dimension of $\zcal$. It corresponds to $\Gr_{m,l}(\C^M\tms\R^L)$. By the set of equations describing $\zcal$, it is obvious that its regular part has the structure of a semiholomorphic foliation.

\subsection{Matrix description}

As it is done with the usual Grassmannian, we describe the CR-Grassmannian as images of linear maps. Consider $A\in \mathsf{M}_{m,l}(\C^M\times\R^L)$ with maximal rank, then
$$
A=\begin{pmatrix}A_1&A_2\\0&A_3\end{pmatrix}
$$
with $A_1$, $A_2$ with complex coefficients, $A_3$ with real coefficients; if $A$ is of maximal rank, then the image of $A$ sits in $\Gr_{m,l}(\C^M\tms\R^L)$.

$A, A'\in \mathsf{M}_{m,l}(\C^M\times\R^L)$ have the same image if there exists an invertible $h\in\mathsf{G}_{m,l}$ such that $A=A'h$.

Incidentally,
$$
\rdim \mathsf{M}_{m,l}(\C^M\times\R^S)=2Mm+2ML+Ll
$$
and
$$
\rdim\mathsf{G}_{m,l}=2m^2+2ml-l^2
$$
so
$$\rdim \mathsf{M}_{m,l}(\C^M\times\R^L)/\mathsf{G}_{m,l}=2Mm+2Ml+Ll-2m^2-2ml-l^2\;.$$

\subsection{Fibration description}

Let $V\in \Gr_{m,l}(\C^M\times\R^L)$, $V_0=V\cap(\C^M\times\{0\})$, $\dim_\C V_0=m$. Fix a linear subspace $S\subset\R^L$ with $\rdim S=l$ and a $\R$-linear injective map $\Phi:S\to\C^M$. Then, let $\Gamma_\Phi=\{(\Phi v,v):v\in S\}$.

We have $V=V_0\oplus \Gamma_\Phi$. Obviously, given $\Phi, \Phi'$,
$$
V_0\oplus\Gamma_\Phi=V_0\oplus\Gamma_{\Phi'}
$$
if and only if $\Gamma_{\Phi'}\subseteq V_0\oplus\Gamma_{\Phi}$, i.e. if $\Phi'(S)\subseteq V_0\oplus \Phi(S)$ (viewing $V_0$ as a subspace of $\C^M$), i.e. if the image of $\Phi-\Phi'$ is contained in $V_0$.

Let
$$
\hcal(E,V_0)=\{\Phi:S\to\C^M\}/\{\Phi:S\to V_0\}
$$
and fix a projection $p:\C^M\times\R^L\to\R^L$, then we have a map
$$
\pi:\Gr_{m,l}(\C^M\times\R^L)\longrightarrow \Gr_m(\C^M)\times \Gr_l(\R^L)
$$
given by $V\mapsto (V\cap \C^M, p(V))$, such that $\pi^{-1}(V_0, S)=\hcal(V_0, S)$.

Again,
$$
\rdim \Gr_m(\C^M)=2m(M-m), \quad \rdim \Gr_l(\R^L)=l(L-l)
$$
and
$$
\rdim \hcal(V_0,S)=2lM-2lm=2l(M-m)\;.
$$

\subsection{The universal CR bundle $U_{m,l}(M,L)$}\label{unibu}
The universal CR bundle ${\sf U}_{m,l}(M,L)$, also denoted ${\sf U}^{M,L}_{m,l}$  is defined by
\be
{\sf U}_{m,l}^{M,L}=\Big\{(V,v)\in \Gr_{m,l}(\C^M\tms\R^L)\tms\C^{M}\tms\R^{L}: v\in V\Big\}.
\ee
Let ${\sf E}\to X$ be a CR bundle of type $(m.l)$ over a CR manifold of type $(n,d)$ and let $f:X\to \Gr_{m,l}(\C^M\tms\R^L)$ denote the continuous map sending the point $x\in X$ to the subspace ${\sf E}_x\sbs\C^{M}\tms\R^{L}$ considered as  an element of $\Gr_{m,l}(\C^M\tms\R^L)$. Then ${\sf E}$ is isomorphic (as a topological bundle) to $f^\ast {\sf U}^{M,L}_{m,l}$. Moreover, if $\{f_t\}_{t\in[0,1]}$ is a homotopy of continuous maps, the bundles $f_t^\ast {\sf U}^{M,L}_{m,l}$ are isomorphic to each other.

If $g:X\to \Gr_{m,l}(\C^M\tms\R^L)$ is a CR map then $g^\ast U^{M,L}_{m,l}$ is a CR bundle of birank $(m,l)$; again, equivalence classes of CR bundles correspond to appropriate homotopy classes of CR maps into $\Gr_{m,l}(\C^M\times\R^L)$.

\medskip

As it is well known, the Oka-Grauert results can be reformulated in terms of homotopies; however, this strategy of proof employs critically the compactness of the classical Grassmannian, whereas the CR Grassmannian is the set of regular points of a (compact) variety, hence non compact. Due to this difficulty, we will move in another direction.

\section{Completeness and vanishing theorems}

\subsection{1-complete foliations}\label{trpsfo}
Let $X$ be a semiholomorphic foliation $X$ of type $(n,d)$ and $\Til X$ its complexification. $X$ is said to be {\it tangentially complete} if $X$ carries a smooth exhaustion function $\phi:X\to\R^+$ which is plurisubharmonic along the leaves and such that $\sup\limits_{X}\phi=+\IN$.

 Let ${\sf N}_{\rm tr}$ be the transverse bundle to the leaves of $X$. A metric on the fibres of $N_{\rm tr}$ is an assignement of a distinguished covering $\{U_j\}$ of $X$ and for every $j$ a smooth map $\l^0_j$ from $U_j$ to the space of symmetric positive $d\tms d$ matrices such that
$$
\l^0_k=\frac{^t\p g_{jk}}{\p t_k}\l^0_j\frac{\p g_{jk}}{\p t_k}.
$$
Denoting $\p$ and $\oli\p$ the complex differentiation along the leaves of $X$, the local tangential forms
\be
2\oli\p\p\log\,\l^0_j-\oli\p\log\,\l^0_j\wedge\p\log\,\l^0_j=\frac{\l^0_j\oli\p\p\l^0_j-2\oli\p\l^0_j\wedge\p\l^0_j}{{\l^0_j}^3}
\ee
\be
\oli\p\p\log\,\l^0_j-\oli\p\log\,\l^0_j\wedge\p\log\,\l^0_j=\frac{2\l^0_j\oli\p\p\l^0_j-3\oli\p\l^0_j\wedge\p\l^0_j}{{\l^0_j}^3}
\ee
actually give global tangential forms $\o$, $\O$.

These forms play the role of the curvature forms of a holomorphic vector bundle on a complex manifold.

In \cites{mt,mt1} the following definitions were given:
\bit
\item $X$ is said to be {\it tranversally 1-complete} ({\it  strongly tranversally 1-complete}) if a metric on the fibres of ${N_{\rm tr}}_{|U}$ can be chosen in such a way that the hermitian form associated to $i\o_{|U}$ ($i\O_{|U}$) has at least $n-q+1$ positive eigenvalues;
\item $X$ is said to be $1\!\!-\!complete$ ({\em strongly} $1\!\!-\!complete$) if it is tangentially and transversally $1$-complete (tangentially and strongly transversally $1$-complete).
\eit
From now on strongly 1-completeness is intended with respect to given exhaustion function $\phi$ and metric on the fibres of $N_{\rm tr}.$

Assume now that $X$ is real analytic, strongly $1\!\!-\!complete$ and let $\Til X$ be the complexification of $X$. Then (see \cite{mt1}*{Theorems 3.1, 3.3}) we have the following

\bit
\item there exist an open neighbourhood $U$ of $X$ in $\Til X$ and a smooth function $u:U\to\R$ such that $u\geq 0$, $u$ is plurisubharmonic in $U$, strongly plurisubharmonic on $U\smi X$ and $X=\{u=0\};$
\item for every $c>0$, $\oli X_c=\{\phi\le c\}$ is a Stein compact of $\Til X$ and every smooth CR function on a neighbourhood of $\oli X_c$ in $X$ can be approximated in the $C^\IN$ topology by smooth CR functions in $X$.
\eit
\subsection{Vanishing theorems}\label{van}
Assume that $\Til X$ is divided by $X$ into two connected components $X_{\pm}$. Let $\oli\ocal_{\pm}$ be the sheaf of germs of holomorphic functions in $X_{\pm}$ smooth up to $\oli X_{\pm}:=X_{\pm}\cup X$ and extend it on $\Til X$ by $0$. Then jump formula for CR functions gives rise to the following (Mayer-Vietoris) exact sequence
\begin{equation}\label{jump}
\xymatrix{0\ar[r]&\ocal\ar[r]&\oli\ocal_+\oplus\oli\ocal_{-}\ar[r]&\scal\ar[r]& 0}
\end{equation}
where $\ocal:=\ocal_{\Til X}$ and $\oli\ocal_+\oplus\oli\ocal_{-}\to\scal$ is defined by $(f\oplus g)=f_{|X}-g_{|X}$.

\nin Then, using the results recalled in \ref{trpsfo} and Kohn's regularity theorem for $\oli\p$ the following was proved in \cite{mt1}*{Theorems 4.1, 4.2}.

\bt If $X$ is real analytic strongly 1-complete semiholomorphic foliation of type $(n,1)$ or a closed, orientable, smooth Levi flat hypersurface in a (connected) Stein manifold $D$ then
$$
H^q(X,\mathcal S)=0
$$
for every $q\geq 1.$\et

\br\label{rem1}
Vanishing for $\scal^\o$ (the sheaf of germs of real analytic CR functions in $X$) fails to be true in general as proved by Andreotti-Nacinovich (see \cite{annac}).
\er
In particular, we obtain the following
\bt\label{teo1}
Let $X$ be a real analytic strongly 1-complete foliation of type ${\rm (}n,1{\rm )}$ and $\mathcal T$ be the sheaf of germs of smooth functions which are locally constant on the leaves. Then
$$
H^q(X,\mathcal T)=0
$$
for $q>n$.
\et
\demo
Let $\O^p_{\rm tg} $ be the sheaf of germs of smooth CR tangential $p$-forms and $\oli\p_{\rm tg}$ the tangential operator $\O^p_{\rm tg}\to\O^{p+1}_{\rm tg}$. In view of the tangential Dolbeault-Grothendieck and \cite{mt1}*{Theorem 5.2}
$$
0\rightarrow \mathcal T \rightarrow
\Omega^0_{\rm tg} \stackrel{\oli\p_{\rm tg}}{\rightarrow}\O^1_{\rm tg}
\stackrel{\oli\p_{\rm tg} }{\rightarrow} \ldots
\stackrel{\oli\p_{\rm tg} }{\rightarrow} \O^n_{\rm tg} \rightarrow 0
$$
is an acyclic resolution of $\mathcal T$. Then, de Rham theorem and $\O^q_{\rm tg}=0$ for $q>n$  imply $H^q(X,\mathcal T)=0$ for $q>n$.
\enddemo
\br\label{rem2}
Cohomology groups $H^q(X,\mathcal T)$ are in general infinite dimensional. They vanish, for instance, if the foliation has a {\em parameters space} i.e. the space of the leaves of $X$ is a smooth curve $T$ (see e.g. \cite{MS}).
\er
\subsubsection{Vanishing for real analytic CR bundles}
The above method still applies to real analytic CR bundles.

 Let $\pi:{\sf E}\to X$ be a real analytic CR bundle of birank ${\rm (}m,l{\rm )}$ and $\ecal$ the sheaf of germs of (smooth)  sections of a ${\sf E}$

 Let $\Til {\sf E}\to \Til X$ be the complexification of ${\sf E}$ and $\Til\ecal$ the sheaf of germs of holomorphic sections of $\Til {\sf E}$. Then it divides $\Til X$ into two domains $\Til X_\pm$ whose boundary is $X$. Denote $\Til\ecal_\pm:=\Til\ecal_{\vert X_\pm}$ and $\oli\ecal_\pm $ the subsheaf of $\Til\ecal_\pm$ on $X_{\pm}$ of germs holomorphic sections which are smooth up to $X$ extended on $\Til X$ by $0$. We have again the following (Mayer-Vietoris) exact sequence
\be
\xymatrix{0\ar[r]&\Til\ecal\ar[r]^\a&\oli\ecal_+\oplus\oli\ecal_-\ar[r]^\b &\Til\ecal_X\ar[r]& 0}
\ee
where $\Til\ecal_X$ is the sheaf of germs of CR section of $\Til{\sf E}_{\vert X}$ and $\b(s\oplus \s)=s_{\vert X}-\s_{\vert X}$.

\bt\label{CO29}
Let $X$ be an orientable analytic strongly 1-complete semiholomorphic foliation of type $(n,1)$, $\pi:{\sf E}\to X$ a real analytic CR bundle of birank ${\rm (}m,l{\rm )}$. Then, if $l=0$
\bit
\item[i)] $H^q(X,\ecal)=0$ for $q\geq1$;
\item[ii)] $\G(X,\ecal)$ generates $\ecal_x$ for every $x\in X$.
\eit
If $l\geq 1$, ${\rm i)}$ and ${\rm ii)}$ are still true provided $X$ satisfies the additional condition
$$
H^1(X,\tcal)=0.
$$
\et
\demo(Sketch) For $l=0$ statement {\rm ii)} is the content of \cite{mt1}*{Theorems 5.2}. In the case $l>0$ we argue as in \cite{mt1}*{Theorems 4.1}.

Let $X=\bigcup\limits_{r\geq1}^{+\IN}X_r$ where $\oli X_r=\{\phi\le r\}$. Every $\oli X_r$ is a Stein compact so, by the methods developed in \cite{gase} and \cite{mt1}*{Theorem 4.1} we show that $H^q(\oli X_r,\Til\ecal_X)=0$ for $q\geq 1$, $r\geq 1$. It follows that $\Til\ecal_{\oli X_r}$ is generated by $\G(\oli X_r,\Til\ecal_{ X_r})$.  Arguing as in \cite{mt1}*{Theorem 4.1} and using Freeman's  approximation theorem (cfr. \cite{FR}*{Theorem 1.3})  finally we obtain the vanishing $H^q(X,\Til\ecal_X)=0$ for $q\geq 1$ and from this it easy to derive the vanishing $H^q(X,\ecal)=0$ for $q\geq 1$. Likewise we prove that $H^q(X,\mcal_x\ecal)=0$ for $q\geq1$ where $\mcal_x$ is the sheaf of germs of CR functions vanishing at a point $x\in\oli X_r$. Then ii) follows as in the classical case.\enddemo
\br
The condition $H^1(X,\tcal)=0$ is necessary. Indeed, let ${\sf E}$ be the trivial bundle $X\tms\R$ considered as CR bundle of type $(0,1)$. Then $\ecal=\tcal$.
\er
\subsubsection{Vanishing for CR bundles}
The general case of a smooth CR bundle $\pi:{\sf E}\to X$ of birank ${\rm (}m,0{\rm )}$ is much more involved. In what follows we just outline the main points of the proof which heavily depends on the thesis work by Sebbar \cites{gase, seb, seb1}.
In order for Sebbar's results to apply in our situation, we examine two coherence conditions.

\medskip

An $\scal_X$-module $\fcal$ is
said to be {\em coherent} if for every point $x$ of $X$  and for
every integer $d \geq 0$ there exist an open neighborhood
$U$ of $x$ and an exact sequence of $\scal_U$-modules
\be
\scal _U^{p_d} \lra \cdots \lra\scal_U^{p_0}
\lra \fcal_U \lra 0,
\ee
where $\scal _U= \scal_X|_U $ and $p_i $  are non-negative integers.

\medskip

An $\oli\ocal_+\oplus\oli\ocal_-$-module $\fcal$ on $\Til X$ is said to be {\em coherent} if for every point $x$ of $X$  and for
every integer $d \geq 0$ there exist an open neighborhood
$U$ of $x$ and an exact sequence of $(\oli\ocal_+\oplus\oli\ocal_-)_U$-modules
\be
(\oli\ocal^{\,\,p_d}_+\oplus\oli\ocal^{\,\,p_d}_ -)_U \lra \cdots \lra(\oli\ocal^{\,\,p_0}_+\oplus\oli\ocal^{\,\,p_0}_-)_U
\lra \fcal_U \lra 0,
\ee
where  $p_i $  are non-negative integers.

\medskip

Let $\fcal$ be a coherent $\scal$-module and extend it by on $\Til X$ by $0$. Then, the jump formula \eqref{jump} for the sheaves $\scal^r$
\be
\xymatrix{0\ar[r]&\ocal^r \ar[r]&{{\oli\ocal}^{\,\,r}_+}\oplus{{\oli\ocal}^{\,\,r}_-}\ar[r]^{\qquad\a}&\scal^r\ar[r]& 0}
\ee
and the results proved by Sebbar in \cite{seb}*{III-3, \S 1} allow us to prove that
\bit
\item[1)] $\fcal$ is a coherent $\oli\ocal_+\oplus\oli\ocal_-$-module;
\item[2)] $H^q(\oli X_c,\fcal)=0$ for every $q\geq1$, $c>0$.
\eit
Arguing as in \cite{mt1}*{Theorems 3.1, 3.3} (using the approximation for CR functions (see \ref{trpsfo})) we then obtain the vanishing
$$
H^q(\oli X,\fcal)=0
$$
for every $q\geq1$.
 In particular
\bt\label{CO291}
Let $X$ be an orientable analytic strongly 1-complete semiholomorphic foliation of type $(n,1)$ and ${\sf E}\to X$ a smooth CR bundle of birank $(m,0)$. Then
\bit
\item[i)] $H^q(X,\ecal)=0$ for $q\geq1$;
\item[ii)] $\G(X,\ecal)$ generates $\ecal_x$ for every $x\in X$.
\eit
\et

\section{The Oka-Grauert principle for real analytic CR bundles}
\bt\label{teo3}
Let $X$ be a strongly 1-complete real analytic semiholomorphic foliation of type ${\rm (}n,1{\rm )}$, ${\sf E}\to X$ a real analytic CR bundle of birank  ${\rm (}m,l{\rm )}$. If ${\sf E}\to X$ is topologically trivial then it is CR trivial.
\et
\demo
Let $\Til X$ be the complexification of $X$ and $\Til {\sf E}$ the complexification of ${\sf E}$ on a neighborod of $X$. If ${\sf E}$ is topologically trivial then $\Til {\sf E}$ is topologically trivial on a neighborod of $X$ (take $n$ continuous sections $s_1,\ldots,s_N:X\to E$ linearly independent at every point $x\in X$ and extend them by continuous sections $\til s_1,\ldots,\til {s}_N$ of $\Til {\sf E}$ on a neighborhood of $X$).

Since $\oli X_c=\{\phi\le c\}$, $c>0$, is a Stein compact of $\Til X$ (\cite{mt1}*{Theorem 3.3}), Grauert theorem implies that $\Til E$ is holomorphically trivial on a neighborhood of $\oli X_c$ in $\Til E$ and so $E$ is CR trivial on a neighborhood of $\oli X_c$ in $X$.

Let $\ecal(U):=\G(U,\ecal)$. Consider $X=\bigcup_{j=1}^{+\IN}X_j$. Then for every $j$ the image of the restriction map $r_j:\ecal(X)\to \ecal(X_j)$ is everywhere dense. For this is enough to show that the restrictions $r^{j+1}_j: \ecal(X_{j+1})\to \ecal(X_j)$ are dense image. Let $s\in\ecal(X_j)$. Since ${\sf E}$ is CR trivial on a neighborhood of $\oli X_{j+1}$ in $X$, there exist CR sections $s_1,s_2,\ldots,s_{m+l}$ on $\oli X_{j+1}$ such that $s=\sum_{k=1}^{m+l}a_ks_k$ where the $a_k^{'s}$ are CR functions on $X_j$. The density of the image of $r_j$ is achieved applying an approximation result \cite{mt1}*{Theorem 3.3}.

To conclude the proof, let $F_j$ be the space of the global sections $s=(s_1,\ldots,s_{m+l})\in\ecal^{\sf cr}(X)^{\oplus(m+l)}$ such that $s_1,\ldots,s_{m+l}$ generate ${\sf E}_{|\oli X_j}$, $j=1,2,\ldots$: clearly all the $F_j$ are open in the Fr\'echet space $\ecal^{\sf cr}(X)^{\oplus(m+l)}$ and $F_1\supset F_2\supset\cdots$; let us prove that they are everywhere dense. Indeed, let $s=(s_1,\ldots,s_{m+l})\in\ecal(X)^{\oplus(m+l)}$ and fix a seminorm $\Vert.\Vert_{\a,K}$ on $\ecal(X)^{\oplus(m+l)}$ with $K=\oli X_{j_0}$, $j_o>j$. Let
$$
U(\a,K,\e)=\Big\{s'\in\ecal(X)^{\oplus{m+l}}:\Vert s'-s\Vert_{\a,K}<\e\Big\}.
$$
Let $\s=(\s_1,\ldots,\s_{m+l})$ where $\s_1,\ldots,\s_{m+l}$ generate $E_{|\oli X_{j_0}}$. Then, $s_p=\sum_{q=1}^{m+l}a_{pq}\s_q$, $1\le p\le m+l$, with $a_{pq}\in\scal(\oli X_{j_0})$. Approximating $\s_1,\ldots,\s_{m+l}$ and the $a_{pq}$ we obtain global sections $\til\s_1,\ldots,\til\s_{m+l}$ which generate ${\sf E}_{|\oli X_{j_0}}$ and such that $\Vert \til\s-s\Vert_{\a,\oli X_{j_0}}<\e$, $\til\s=(\til\s_1,\ldots,\til\s_{m+l})$, and $\til\s\in\ecal^{\sf cr}(X)^{\oplus{m+l}}$ i.e. $\til\s\sbs \Big(F_{j_0}\cap U(\a,K,\e)\Big)\sbs\Big(F_j\cap U(\a,K,\e)\Big)$. This show that sets $F_j$, $j\geq 1$ are everywhere sense in $\ecal(X)^{\oplus{m+l}}$, so, in view of Baire Lemma, $\cap_{j\geq1}F_j\neq\ES$. It follows that there are global sections $\s_1,\ldots,\s_{m+l}$ which generate $\ecal(X)^{\oplus{m+l}}$ and consequently that ${\sf E}$ is CR trivial.
\enddemo
\bt\label{teo31}
Let $X$ be a strongly 1-complete real analytic semiholomorphic foliation of type $\(n,d\)$. Then two topologically equivalent real analytic CR bundles of birank $\(m,0\)$ are CR equivalent.
\et
\demo
Let ${\sf E}\to X$, ${\sf F}\to X$ be topologically equivalent CR bundles of birank $\(m,0\)$. Let $\ucal=\{U_j\}_{j\in\N}$ be covering of $X$ by balls of $X$ such that ${\sf E}_{|U_j}$ and ${\sf F}_{|U_j}$ are trivial for every $j\in\N$ and $\{h_{ij}:U_i\cap U_j\to\C^\ast\}$, $\{g_{ij}:U_i\cap U_j\to\C^\ast\}$ the cocycles of ${\sf E}$, ${\sf F}$ respectively, associated to $\ucal$.

Let $\Til {\sf E}\to\Til X$, $\Til {\sf F}\to\Til X$ be the respective complexifications; then, they are topologically equivalent on a neighbourhood of $X$ in $\Til X$.


Let $X=\bigcup\limits_{r\geq1}^{+\IN}X_r$ where $\oli X_r=\{\phi\le r\}$. Since $\oli X_r$ is a Stein compact (see subsection \ref{van}), by Grauert Theorem for every $r$ there exists a Stein neighborhood $U_r$ of $\oli X_r$ in $\Til X$ such that $\Til {\sf E}_{\vert U_r}\stackrel{\rm hol}{\simeq}\Til {\sf F}_{\vert U_r}$. Moreover, for every compact subset $K\sbs X$ there exist sections $u_1,u_2,\ldots,u_N\in\ecal_{\vert K}$, $N=N(K)$, which generate $\ecal$ on $K$.

Consider the sheaf ${\sf Hom}(\ecal,\fcal)$; it is a Fr\'echet sheaf. Every local  isomorphism ${\sf E}\stackrel{\rm cr}{\simeq}{\sf F}$ defines a local section of ${\sf Hom}(\ecal,\fcal)$.

Let $s$ be a section in ${\sf Hom}(\ecal,\fcal)(U)$, $U\sbs X$ which defines an isomorphism $\ecal_{\vert U}\simeq\fcal_{\vert U}$. If $s'$ near $s$ (in the Fr\'echet topology on the space of sections ${\sf Hom}(\ecal,\fcal)(U)$) then $s'$ is a local isomorphism.

Let $s_1\in{\sf Hom}(\ecal,\fcal)(\oli X_1)$ be the section determined by the isomorphism $\Til {\sf E}_{\vert U_1}\stackrel{\rm hol}{\simeq}\Til {\sf F }_{\vert U_1}$. Then there exists $\delta_1$ such that if a section $s'\in{\sf Hom}(\ecal,\fcal)(\oli X_1)$ satisfies $\Vert s'-s_1\Vert_{\oli X_1}<\delta_1$, where $\Vert s'-s_1\Vert$ is the distance in the Fr\'echet topology; then $s'$ is an isomorphism on $\oli X_1$.

Fix sections $u_1,u_2,\ldots,u_N\in\ecal_{\vert \oli X_2 }$, which generate $\ecal$ on $\oli X_2$ and let $s'\in{\sf Hom}(\ecal,\fcal)(\oli X_2)$ be an isomorphism $\ecal_{\vert\oli X_2}{\simeq}\fcal_{\vert\oli X_2}$ (determined by the isomorphism $\Til {\sf E}_{\vert U_2}\stackrel{\rm hol}{\simeq}\Til {\sf F}_{\vert U_2}$). Then $s'(u_1),s'(u_2),\ldots,s'(u_m)$ generate $\fcal_{\vert\oli X_2}$ hence on one has
$$
s(u_i)=\sum_{k=1}^mA'_{ki}s'(u_k)
$$
where the $A'_{ki}$ are (smooth) CR functions on $\oli X_1$ and the matrix $A':=(A'_{ki})$ is of maximal rank $m$ on $\oli X_1$.

Since $\oli X_2$ is $\ocal (U_2)$-convex we can approximate $A'$ by matrices of holomorphic functions of rank $m$ on $\oli X_1$. In view of\cite{forst}*{Theorem 5.4.4}, $A'$ approximates by matrices of holomorphic functions of rank $m$ on $\oli X_2$ so that we can choose a matrix $A$ of   CR functions of maximal rank $m$ on $\oli X_2$ in such a way to have $\Vert As'-s_1\Vert<2^{-1}\delta_1$. Set $s_2:=As'$, then $s_2$ is an isomorphism on $X_2$. Let $\delta_2<\delta_1$  and $s_3$ be an isomorphism on $X_3$ such that  $\Vert s_3-s_2\Vert<2^{-1}\delta_2$ and so on. By this procedure for every $k\in\N$ we construct a section $s_k\in{\sf Hom}(\ecal,\fcal)(\oli X_k)$ which is a isomorphism on $\oli X_k$ and
$$
\Vert s_{k+1}-s_k\Vert<2^{-k}\delta_k
$$
and $\delta_1>\delta_2>\cdots$.
%
By standard techniques, $s_k\to s\in{\sf Hom}(\ecal, \fcal)(X)$ and for every $k$ the restriction of $s$ to $\oli X_k$ is an isomorphism. This shows that $s$ is an isomorphism $\ecal\stackrel{\rm CR}{\simeq}\fcal$.
\enddemo

With the techniques developed in the next section, we can proof an analogue result for $(m,1)$-bundles, with an extra assumption:
\bt\label{teo32} Let $X$ be a strongly 1-complete real analytic semiholomorphic foliation of type $\(n,d\)$.
If
$$
H^1(X,\tcal)=H^1(X,\Z)=0\;,
$$
then two topologically equivalent real analytic CR bundles of birank $(m,1)$ are CR equivalent.
\et


\section{The Oka-Grauert principle for smooth CR bundles}\label{O-Gra2}

\subsection{The case of birank $(1,0)$-bundles}\label{O-Gra1}
We have the following
\begin{theorem}\label{teo1OG}
Let $X$ be a strongly 1-complete real analytic semiholomorphic foliation of type $\(n,1\)$. Then
$$
{\sf Vect}_{\,\sf cr}^{(1,0)}(X)\simeq
{\sf Vect}_{\,\sf top}^{(1,0)}(X).
$$
\end{theorem}
\demo
The proof is the same as that of Oka for Stein manifolds.  We know that
$$
{\sf Vect}_{\,\sf top}^{(1,0)}(X)\simeq H^1(X,\ccal^\ast), \>\>{\sf Vect}_{\,\sf cr}^{(1,0)}(X)\simeq H^1(X,\scal^\ast)
$$
where $\ccal=\ccal_X$ is the sheaf of germs of continuous functions in $X$).

Consider the exact diagram of sheaves
$$
\xymatrix{0\ar[r]&\Z\ar[r]^{}\ar[d]_{} &\scal\ar[r]^{\exp 2\pi i}\ar[d]_{} & \scal^\ast\ar[r]\ar[d]_{}&0\\0\ar[r]&\Z\ar[r]^{}& \ccal\ar[r]^{\exp 2\pi i}&\ccal^\ast\ar[r]&0}
$$
where $\ccal=\ccal_X$ is the sheaf of germs of continuous functions in $X$. Since $H^q(X,\ccal)=0$, $q\geq 1$, and $H^q(X,\scal)=0$, $q\geq 1$, (see \cite{mt1}*{Theorems 4.1}) passing to cohomology we get the following exact diagram
of groups
$$
\xymatrix{0\ar[r]^{}\ar[d]_{}&H^1(X,\scal^\ast)\ar[r]^{}\ar[d]_{} &H^2(X,\Z)\ar[r]^{}\ar[d]_{\simeq} & 0\ar[d]_{} \\0\ar[r]&H^1(X,\ccal^\ast)\ar[r]^{} & H^2(X,\Z)\ar[r]^{} & 0}.
$$
whence
$$
{\sf Vect}_{\,\sf cr}^{(1,0)}(X)\simeq H^1(X,\scal^\ast)\simeq H^1(X,\ccal^\ast)\simeq {\sf Vect}_{\,\sf top}^{(1,0)}(X).
$$
In particular
\bit
\item[a)] two topologically equivalent CR vector bundles of birank $\(1,0\)$ over $X$ are CR equivalent;
\item[b)] every topological vector bundle of type $(1,0) $ over $X$ admits an equivalent CR vector bundle structure.
\eit
\enddemo

\subsection{The case of $(1,1)$-bundles}\label{ssc11}
The result we aim at is the following.
\begin{theorem}\label{teo2}
Let $X$ be a strongly 1-complete real analytic semiholomorphic foliation of type $\(n,1\)$. Assume that
$$
H^1(X,\tcal)=H^1(X,\Z)=0.
$$
Then
$$
{\sf\epsilon}_{\!_X}:{\sf Vect}_{\,\sf cr}^{(1,1)}(X)\longrightarrow
{\sf Vect}_{\,\sf top}^{(1,1)}(X)
$$
is a bijection.
\end{theorem}

However, we will not present a complete proof, but only a reduction to a cohomological problem, with the same general strategy presented by Cartan in \cite{Cart} for the classical Oka-Grauert results. To be more precise,
\begin{enumerate}\renewcommand{\theenumi}{\roman{enumi}}
\item in Proposition \ref{CRtrivial} we show that a topologically trivial CR bundle of birank $(m,1)$ is CR trivial,
\item in Proposition \ref{CRinj} we show that the map $ {\sf\epsilon}_{\!X}$ is injective,
\item in the last part of subsection \ref{ssc11}, we reduce the surjectivity to a cohomological problem.
\end{enumerate}\renewcommand{\theenumi}{\arabic{enumi}}
We do not present the proof of the cohomological property needed to conclude, but we believe it can be achieved by techniques similar to the ones employed in the classical case in \cite{Cart}.

\bigskip

We first prove the following
\bl\label{rid}
Let ${\sf E}\to X$, ${\sf F}\to X$ topologically equivalent CR bundles of birank $\(m,1\)$. Then they can be defined by cocycles $\{h_{ij}\}$, $\{g_{ij}\}$ of the form
\begin{equation}\label{CO22}
h_{ij}=\left( \begin{array}{cc}
	F_{ij}&B_{ij}\\
	0&1
\end{array} \right),\>\>
g_{ij}=\left( \begin{array}{cc}
	F_{ij}&G_{ij}\\
	0&H_{ij}\end{array} \right).
\end{equation}
\el
\demo
Suppose that ${\sf E}\to X$, ${\sf F}\to X$ are topologically equivalent CR bundles of birank $\(m,1\)$ and let
\begin{equation}\label{CO2}
h_{ij}=\left( \begin{array}{cc}
	A_{ij}&B_{ij}\\
	0&C_{ij}
\end{array} \right),\>\>
g_{ij}=\left( \begin{array}{cc}
	F_{ij}&G_{ij}\\
	0&H_{ij}
\end{array} \right)
\end{equation}
the cocycles of ${\sf E}$ and ${\sf F}$ respectively.

By hypothesis there are continuous maps $a_i;U_i\to{\sf G}_{1,1}$
\be
a_i=\left( \begin{array}{cc}
	M_i&S_i\\
	0&P_i
\end{array} \right)
\ee
such that
\be
h_{ij}=a^{-1}_ig_{ij}a_j
\ee
on $U_{ij};=U_i\cap U_j$ where $P_i$ is constant on the leaves.

In particular, (\ref{CO2}) gives the following identities
\begin{equation}
\begin{cases}\label{CO4}
A_{ij}=M_i^{-1}F_{ij}M_j\\
C_{ij}=P_i^{-1}H_{ij}P_j.
\end{cases}
\end{equation}
Then, again by \cite{mt1}*{Theorems 4.1, 4.2}, we can suppose that $M_i$ and $M_j$ are CR functions.
So, $\{h_{ij}\}$ is CR equivalent to the cocycle
\be
\left( \begin{array}{cc}
	M_i&0\\
	0&1
\end{array} \right)\left( \begin{array}{cc}
	A_{ij}&B_{ij}\\
	0&C_{ij}\end{array} \right)\left( \begin{array}{cc}
	M_j^{-1}&0\\
	0&1
\end{array} \right)=\left( \begin{array}{cc}
	F_{ij}&\star\\
	0&C_{ij}
\end{array} \right)
\ee
thus we may assume ${\sf E}$ and ${\sf F}$ are defined by the CR cocycles
\be
h_{ij}=\left( \begin{array}{cc}
	F_{ij}&B_{ij}\\
	0&C_{ij}
\end{array} \right),\>\>
g_{ij}=\left( \begin{array}{cc}
	F_{ij}&G_{ij}\\
	0&H_{ij}
\end{array} \right).
\ee
respectively. From here it is easy to conclude. \enddemo

By Lemma \ref{rid} we may suppose that ${\sf E}$, ${\sf F}$ are defined by the cocycles \eqref{CO22}.

If  $\{h_{ij}\}\stackrel{\rm top}{\simeq}\{g_{ij}\}$,$\{H_{ij}\}$ is coboundary with values in $\R^\ast$ so $\{H_{ij}\}=\exp(H_i-H_j)$ with $H_i:U_i\to\C$,  $H_j:U_j\to\C$ continuous and constant on the leaves. Since $F_{ij}/H{ij}$ is real valued function
\be
{\sf Im}(H_i-H_j)=\pi N_{ij}
\ee
with $N_{ij}\in\Z$. Clearly, $\{N_{ij}\}$ is a 1-cocycle with values in the constant sheaf $\Z$ hence $N_{ij}=N_i-N_j$, $N_i,N_j\in\Z$, by the hypothesis $H^1(X,\Z)=0$. Then ${\sf Im}\,H_i-\pi N_i={\sf Im}\,H_j-\pi N_j$ for evey $i,j$, so $\{{\sf Im}\,H_i-\pi N_i\}$ is a continuous function $u$ in $X$ and

\beqnn
H_{ij}&=&\exp({\sf Re}\,H_i+i\pi N_i+iu)\exp(-({\sf Re}\,c_j+i\pi N_j+iu)=\\
&&\exp({\sf Re}\,H_i+i\pi N_i))\exp(-({\sf Re}\,H_j+i\pi N_j).
\nonumber\\
\eeqnn
It follows that
\be
\exp(i(N_j-N_i))H_{ij}=\exp H_i\exp(-H_j)
\ee
whence that the function $H_i-H_j$ is smooth and constant on the leaves. Clearly. $\{H_i-H_j\}$ is a 1-cocycle with values in the sheaf $\tcal$ so, by hypothesis $H^1(X,\tcal)=0$, there exist smooth functions $\psi_i:U_i\to\R$ which are constant on the leaves and such that $H_i-H_j=\psi_i-\psi_j$. Then the smooth functions $P_i:=\exp(\psi_i+iN_i\pi)$ are real valued, constant on the leaves and
\be
H_{ij}=\frac{\exp(\psi_i+iN_i\pi)}{\exp(\psi_j+iN_j\pi)}=\frac{P_i}{P_j}.
\ee
It follows that
\be
g_{ij}=\left( \begin{array}{cc}
	F_{ij} &G_{ij}\\
	0&P_i^{-1}/P_j
\end{array} \right)=\left( \begin{array}{cc}
	1&0\\
	0&P_i^{-1}\end{array}\right)\left( \begin{array}{cc}
	F_{ij}&G_{ij}\\
	0&1\end{array}\right)\left( \begin{array}{cc}
	1&0\\
	0&P_j\end{array}\right).
\ee
Therefore we may assume that
\begin{equation}\label{CO221}
h_{ij}=\left( \begin{array}{cc}
	F_{ij}&B_{ij}\\
	0&1
\end{array} \right),\>\>
g_{ij}=\left( \begin{array}{cc}
	F_{ij}&G_{ij}\\
	0&1\end{array} \right).
\end{equation}

\bp \label{CRtrivial}Let ${\sf E}\to X$ be a topologically trivial CR bundle of birank $\(m,1\)$.  Assume that
$$
H^1(X,\tcal)=H^1(X,\Z)=0.
$$
Then ${\sf E}\to X$  is CR trivial.\ep
\demo With the notation introduced above, if $\{g_{ij}\}$ is topologically trivial then $F_{ij}$ is also topologically trivial, and consequently, since $H^{-1}(X,\scal)=0$, we have  $F_{ij}=F_i^{-1}F_j$ for some CR functions  $F_i$. It follows that
\be
\left( \begin{array}{cc}
	F_i^{-1}&0\\
	0&1\end{array}\right)\left( \begin{array}{cc}
	F_{ij}&B_{ij}\\
	0&1\end{array}\right)\left( \begin{array}{cc}
	F_j&0\\
	0&1\end{array}\right)=\left( \begin{array}{cc}
	1&B'_{ij}\\
	0&1
\end{array} \right):=h'_{ij}
\ee
i.e. $\{h'_{ij}\}\stackrel{\rm cr}{\simeq}\{h_{ij}\}$. Cocycle condition for $h'_{ij}$ implies that $B'_{ij}$ is a 1-cocycle with values in $\scal$ wence a coboundary: $B'_{ij}=B_i-B_j$. It follows that
$$
 \left( \begin{array}{cc}
	1&B_i\\
	0&1
\end{array} \right)^{-1}
\left( \begin{array}{cc}
	1&B'_{ij}\\
	0&1
\end{array} \right)
\left( \begin{array}{cc}
	1&B_j\\
	0&1
\end{array} \right)=
\left( \begin{array}{cc}
1&0\\
	0&1
\end{array} \right)
$$
that ${\sf E}$ is CR trivial.\enddemo

Let us prove now that $\{h_{ij}\}\stackrel{\rm top}{\simeq}\{g_{ij}\}$ implies $\{h_{ij}\}\stackrel{\rm cr}{\simeq}\{g_{ij}\}$. We can assume that the cocycles are given by \eqref{CO221}.

By hypothesis there are continuous maps $c_i:U_i\to{\sf G}_{1,1}$
\be
c_i=\left( \begin{array}{cc}
	F_i&G_i\\
	0&H_i
\end{array} \right)
\ee
such that
\be
h_{ij}=c_i^{-1}g_{ij}c_j
\ee
on $U_i\cap U_j$ where $H_i$ is constant on the leaves.  Moreover, in view of the particular form of $h_{ij}$ and $g_{ij}$, we have $F_i=F_j=F$, $H_i=H_j=H$, where $F$ and $H$ are continuous  global functions in $X$ and $H$ is an invertible real valued function constant on the leaves. Moreover, multiplying by $H^{-1}$ we may assume
\begin{equation}\label{CO241}
c_i=\left( \begin{array}{cc}
	F&G_i\\
	0&1
\end{array} \right);
\end{equation}
thus
\begin{equation}\label{CO242}
c_i^{-1}c_j=\left( \begin{array}{cc}
	1&-F^{-1}G_i+F^{-1}G_j\\
	0&1
\end{array} \right)
\end{equation}
From \eqref{CO22} and \eqref{CO241}, we deduce that $\sf E$ and $\sf F$ have structure group $\sf{G}:=\C^\ast\tms\C$, identified with the group of the matrices
\be
\left( \begin{array}{cc}
	a&b\\
	0&1
\end{array} \right)
\ee
$a\in\C^{\ast}$, $b\in\C$ with product $(a,b)\cdot (c,d)=(ac,ad+b)$

The Lie algebra $\mathfrak{g}$ of $\sf{G}$ is $\C\tms\C$ with bracket
$$
[(\z,\eta),(\z',\eta')]=(0,\z\eta'-\eta\z').
$$
If $\xi=(\z,\eta)\in \gfra$, $\l\in\C$, then $\exp \l\xi=(\exp \l\z,\l\eta).$ Then $\{c_i\}$ gives a section of the CR bundle $Z\to X$  and we have to prove that
\begin{quote}
(ii) every continuous section $X\to Z$ is homotopic to a CR section.
\end{quote}
This reduces to a cohomological problem. 

We follow \cite{Cart}*{\bf 4} taking into account \eqref{CO242}. Let ${\sf G} _0=\{1\}\tms\C\tms\{1\}\simeq\C$ ($(1,a,1)\cdot(1,b,1)=(1,a+b,1)$, $(1,a,1)^{-1}=(1,-a,1)$ in ${\sf G} _0$). For every open subset $U$ of $X$, let $\fcal(U)$ be the topological group of the continuous maps $f:I\to C^0(U,{\sf G} _0)$ satisfying: $f(0)={\sf 1}$ and $f(1)\in C^{{\sf cr}}(U,{\sf G}_0)$.
Both are topological groups (with a Fréchet topology).

Then, arguing as in \cite{Cart}*{\bf 4}, (ii) will be a consequence of the following
\bt\label{teo221}
If $X$ is a real analytic strongly 1-complete semiholomorphic foliation of type $(n,1)$. Then
$$
H^1(X,\fcal)=0.
$$
\et
\demo
We have to prove the following: given an open covering $\{U_i\}_i$ of $X$ and continuous functions $a_{ij}:U_{ij}\tms I:\to\C$, $a_{ij}=a_{ij}(x,\l)$, such that
\be
a_{ij}+a_{jk}+a_{ki}=0,\>\>, a_{ij}(\cdot,0)=0, \>\>a_{ij}(\cdot,1)\in\scal(U_{ij})
\ee
then there exist continuous functions $a_i:U_i\tms I\to \C$ such that
\be
a_{ij}=a_i-a_j,\>\>, a_i(\cdot,0)=0, \>\>a_i(\cdot,1)\in\scal(U_i).
\ee
The proof of this theorem is achieved adapting to our case the method followed in \cite{Cart}*{\bf 4}, using the results recalled in Section \ref{trpsfo}.\enddemo

Thus, we obtained the following
\bp \label{CRinj} In our hypotheses, the map
$$
{\sf\epsilon}_{\!_X}:{\sf Vect}_{\,\sf cr}^{(1,1)}(X)\longrightarrow
{\sf Vect}_{\,\sf top}^{(1,1)}(X)
$$
is injective.\ep

\bigskip

It remains to prove that ${\sf\epsilon}_{\!_X}$ is onto i.e. if $\{h_{ij}\}$ is a continuous cocycle with values in $\sf G_{1,1}$, $h_{ij}:U{ij}\to \sf G_{1,1}$ then there exist continuous maps $c_i:U_i\to\sf G_{1,1}$ such that the cocycle $c_i^{-1}h_{ij}c_j$ is CR.

We may assume that
 \be
h_{ij}=\left( \begin{array}{cc}
	F_{ij}&B_{ij}\\
	0&1
\end{array} \right).
\ee
Let $\sf G=\C^\ast\tms\C\tms\{1\}\simeq\C^\ast\tms\C$ (with product $(a,b)\cdot(a',b')=(aa',ab'+b)$ and $(a,b)^{-1}=(a^{-1},-a^{-1}b)$). Denote $\gcal^{\rm cont}$ (respectively $\gcal^{\rm cr}$) the sheaf of germs of continuos (respectively CR) maps with values in $\sf G$. Then showing that the cocycle $c_i^{-1}h_{ij}c_j$ is CR is equivalent to prove that the natural map
\begin{equation}\label{eq_teo2}
H^1(X,\gcal^{\rm cr})\longrightarrow H^1(X,\gcal^{\rm cont})
\end{equation}
is onto.
\bigskip

So, we have reduced the proof of Theorem \ref{teo2} to proving that the map \eqref{eq_teo2} is surjective.

\subsection{The general case} The case $(m,l)$ essentially reduces to the cases $(m,0)$, $(0,l)$. Indeed, let $\sf E,\sf F\in{\sf Vect}_{\,\sf cr}^{(m,l)}(X)$ be topologically equivalent with respective   cocycles
\be
h_{ij}=\left( \begin{array}{cc}
	A_{ij}&B_{ij}\\
	0&C_{ij}
\end{array} \right),\>\>
g_{ij}=\left( \begin{array}{cc}
	F_{ij}&G_{ij}\\
	0&H_{ij}\end{array} \right).
\ee
Then $\{A_{ij}\}$, $\{F_{ij}\}$ are topologically equivalent CR cocycles of type $(m,0)$ and $\{C_{ij}\}$, $\{H_{ij}\}$ of type $(0,l)$ (see \eqref{CO4})
%
%
%
Moreover, $\{A_{ij}\}$, $\{F_{ij}\}$ are are 1-cocycle with values in ${\rm GL(m,\C)}$ with, by hypothesis, are topologically equivalent. If the case $(m,0)$ is already solved there exist CR maps $M_i, P_i$ such that
$$
F_{ij}=M_i^{-1}A_{ij}M_j,\>\>H_{ij}=P_i^{-1}C_{ij}P_j
$$
It follows that $\{g_{ij}\}$ is CR equivalent to
Taking
\be
c^{-1}_i=\left( \begin{array}{cc}
	I&0\\
	0&P^{-1}_i\end{array} \right)
\left( \begin{array}{cc}
	M^{-1}_i&0\\
	0&I\end{array} \right),\>\>c_j=\left( \begin{array}{cc}
	M_j&0\\
	0&I\end{array} \right)\left( \begin{array}{cc}
	I&0\\
	0&P_j\end{array} \right)
\ee
we have
\be
c^{-1}g_{ij}c_j=\left( \begin{array}{cc}
	A_{ij}&\star\\
	0&C_{ij}\end{array} \right)
	\ee
so we may assume
\begin{equation*}
h_{ij}=\left( \begin{array}{cc}
	A_{ij}&B_{ij}\\
	0&C_{ij}
\end{array} \right),\>\>
g_{ij}=\left( \begin{array}{cc}
	A_{ij}&G_{ij}\\
	0&C_{ij}\end{array} \right).
\end{equation*}
Now, topological equivalence implies that there are continuous maps $c_i:X\to{\sf G}_{m,l}$ \begin{equation*}
c_i=\left( \begin{array}{cc}
	F&G_i\\
	0&H
\end{array} \right)
\end{equation*}
such that
\be
h_{ij}=c^{-1}_ig_{ij}c_j.
\ee
Moreover
\begin{equation}\label{CO227}
c^{-1}_ic_j=\left( \begin{array}{cc}
	I&-F^{-1}G_i+F^{-1}G_j\\
	0&I
\end{array} \right).
\end{equation}
$c=\{c_i\}$ gives a section of the CR bundle $Z\to X$ and we have to prove that
\begin{quote}
(ii) $c$ is homotopic to a CR section.
\end{quote}
This reduces again to a cohomological problem. 

We follow \cite{Cart}*{\bf 4} taking into account \eqref{CO227}. Let ${\sf G} _0=\{{\bf 1}\}\tms\C^{m+l}\tms\{{\bf 1}\}\simeq\C^{m+l}$ ($({\bf 1},a,{\bf1})\cdot({\bf1},b,{\bf 1})=({\bf 1},a+b,{\bf 1})$, $({\bf 1},a,{\bf 1})^{-1}=({\bf 1},-a,{\bf 1})$ in ${\sf G} _0$). Let ${\rm I}=[0,1]$. For every open subset $U$ of $X$, let $\fcal(U)$ be the topological group of the continuous maps $f:{\rm I}\to C^0(U,{\sf G} _0)$ satisfying: $f(0)={\sf I}$ and $f(1)\in C^{{\sf cr}}(U,{\sf G}_0)$.
Both are topological groups (with a Fr$\acute{e}$chet topology).

Injectivity of the map
$$
{\sf\epsilon}_{\!_X}:{\sf Vect}_{\,\sf cr}^{(m,l)}(X)\longrightarrow
{\sf Vect}_{\,\sf top}^{(m,l)}(X)
$$
will be a consequence of the following
\bt\label{teo2211}
If $X$ is a real analytic strongly 1-complete semiholomorphic foliation of type $(n,1)$. Then
$$
H^1(X,\fcal)=0.
$$
\et
As for surjectivity denote $\gcal^{\rm cont}$ respectively $\gcal^{\rm cr}$  the sheaf on $X$ of germs of continuous maps with values in ${\sf G}_{m,l}$ and that one of germs of CR maps with values in ${\sf G}_{m,l}$.

Then it must be prove that
$$
H^1(X,\gcal^{\rm cr})\lra H^1(X,\gcal^{\rm cont})
$$
is onto.

Observe that if the cases $\(m,0\)$, $\(0,l\)$ are solved we may assume that every cohomology class in $H^1(X,\gcal^{\rm cont})$ is represented by a cocycle
\be
h_{ij}=\left( \begin{array}{cc}
	A_{ij}&B_{ij}\\
	0&C_{ij}
\end{array} \right)
\ee
where $A_{ij}$ and $C_{ij}$ are CR, $B_{ij}$  is continuous and satisfies the cocycle condition
\be
A_{ij}B_{ji}+
	B_{ij}C_{ji}=0.
\ee


\begin{bibdiv}
\begin{biblist}
%

   \bib{annac}{book}{
     author={Andreotti, A.},
   author={Nacinovich, M.},
   title={Analytic convexity and the principle of Phragm\'en-Lindel\"of},
   note={Pubblicazioni della Classe di Scienze: Quaderni. [Publications of
   the Science Department: Monographs]},
   publisher={Scuola Normale Superiore Pisa, Pisa},
   date={1980},
   pages={182},
       }
   \bib{Cart}{article}{
    author={Cartan, Henri},
   title={Espaces fibr\'es analytiques},
   language={French},
   conference={
      title={Symposium internacional de topolog\'\i a algebraica International
      symposium on algebraic topology},
   },
   book={
      publisher={Universidad Nacional Aut\'onoma de M\'exico and UNESCO, Mexico
   City},
   },
   date={1958},
   issn={},
   review={},

 }

%
   \bib{forst}{book}{
    author={Forstneri\v c, Franc},
   title={Stein manifolds and holomorphic mappings},
   series={Ergebnisse der Mathematik und ihrer Grenzgebiete. 3. Folge. A
   Series of Modern Surveys in Mathematics [Results in Mathematics and
   Related Areas. 3rd Series. A Series of Modern Surveys in Mathematics]},
   volume={56},
   note={The homotopy principle in complex analysis},
   publisher={Springer, Heidelberg},
   date={2011},
   pages={xii+489},
   issn={},
   review={},
 }

\bib{FR}{article}{
    author={Freeman, Michael},
   title={Tangential Cauchy-Riemann equations and uniform approximation},
   journal={Pacific J. Math.},
   volume={33},
   date={1970},
   pages={101--108},
   issn={},
   review={},
   }
   \bib{gase}{article}{
   author={Gay, Roger},
   author={Sebbar, Ahmed},
   title={Division et extension dans l'alg\`ebre $A^\infty(\Omega)$ d'un
   ouvert pseudo-convexe \`a bord lisse de ${\bf C}^n$},
   language={French},
   journal={Math. Z.},
   volume={189},
   date={1985},
   number={3},
   pages={421--447},
   issn={},
   review={},
  }

\bib{G1}{article}{
   author={Grauert, Hans},
   title={Holomorphe Funktionen mit Werten in komplexen Lieschen Gruppen},
   language={German},
   journal={Math. Ann.},
   volume={133},
   date={1957},
   pages={450--472},
}
\bib{G2}{article}{
   author={Grauert, Hans},
   title={Analytische Faserungen \"uber holomorph-vollst\"andigen R\"aumen},
   language={German},
   journal={Math. Ann.},
   volume={135},
   date={1958},
   pages={263--273},
}
\bib{G3}{article}{
   author={Grauert, Hans},
   title={Approximationss\"atze f\"ur holomorphe Funktionen mit Werten in
   komplexen R\"amen},
   language={German},
   journal={Math. Ann.},
   volume={133},
   date={1957},
   pages={139--159},
}
	\bib{gro}{article}{
   author={Gromov, M.},
   title={Oka's principle for holomorphic sections of elliptic bundles},
   journal={J. Amer. Math. Soc.},
   volume={2},
   date={1989},
   number={4},
   pages={851--897},
}

%
   \bib{MS}{book}{
   author={Moore, Calvin C.},
   author={Schochet, Claude L.},
   title={Global analysis on foliated spaces},
   series={Mathematical Sciences Research Institute Publications},
   volume={9},
   edition={2},
   publisher={Cambridge University Press, New York},
   date={2006},
   pages={xiv+293},
   isbn={978-0-521-61305-7},
   isbn={0-521-61305-1},
   review={\MR{2202625}},
} \bib{mt1}{article}{
   author={Mongodi, Samuele},
   author={Tomassini, Giuseppe},
   title={1-complete semiholomorphic foliations},
   journal={Trans. Amer. Math. Soc.},
   volume={368},
   date={2016},
   number={9},
   pages={6271--6292},
   issn={},
   review={},

   }
   \bib{mt}{article}{
  author={Mongodi, Samuele},
   author={Tomassini, Giuseppe},
   title={Transversally pseudoconvex semiholomorphic foliations},
   journal={Atti Accad. Naz. Lincei Rend. Lincei Mat. Appl.},
   volume={26},
   date={2015},
   number={1},
   pages={23--36},
   issn={},
   review={},
   }

%
   \bib{oka}{article}{
    Author = {Oka, Kiyoshi},
    Title = {{Sur les fonctions analytiques de plusieurs variables. III. Deuxi\`eme probl\`eme de Cousin.}},
    Journal = {{J. Sci. Hiroshima Univ., Ser. A}},
    Volume = {9},
    Pages = {7--19},
    date = {1939},
    Publisher = {Hirosima University, Hiroshima},
    Language = {French},
}
   \bib{rea}{article}{
   author={Rea, C.},
   title={Levi-flat submanifolds and holomorphic extension of foliations},
   journal={Ann. Scuola Norm. Sup. Pisa (3)},
   volume={26},
   date={1972},
   pages={665--681},
   issn={},
   review={},
   }
    \bib{seb}{thesis}{
   author={Sebbar, A.},
   author={},
   title={Espaces fibrés $\acal^\IN$ et Théorème de Grauert},
   language={},
   type={Thèse presentée à l'Université de Grenoble},
   volume={},
   date={1985},
   language={French},
   pages={},
   issn={},
   review={},
}
   \bib{seb1}{article}{
  author={Sebbar, A.},
   title={Principe d'Oka-Grauert dans $A^\infty$},
   language={French},
   journal={Math. Z.},
   volume={201},
   date={1989},
   number={4},
   pages={561--581},
   issn={},
   review={},
}

\end{biblist}
\end{bibdiv}
\end{document}